\documentclass[11pt]{article}
\makeatletter \@addtoreset{figure}{section} \makeatother
\makeatletter
\long\def\@makecaption#1#2{%
   \vskip 10\p@
   \setbox\@tempboxa\hbox{{#1}\ \ #2}%
   \ifdim \wd\@tempboxa >\hsize
       {#1}\ \ #2\par
   \else
       \hbox to\hsize{\hfil\box\@tempboxa\hfil}%
   \fi}
\makeatother \makeatletter
\newcommand
\myroman[1]{\@roman #1}
\makeatother

\newtheorem{thm}{Theorem}[section]

\newtheorem{lem}[thm]{Lemma}

\newcommand{\qed}{{\hfill\rule{3pt}{7pt}}}
\def\pf{\noindent {\it Proof.} }

\def\qed{\hfill \rule{4pt}{7pt}}
\def\pf{\noindent {\it Proof} }
\title{\bf Existence of Nontrivial Solutions for   $p$-Laplacian Equations in $\mathbf{R}^{N}$}
\author
{
Chungen Liu
\footnote
{School of Mathematical Sciences and LPMC,  Nankai University,  Tianjin 300071,  China.\newline\indent \ \
Partially  supported by NFS of China
and 973 Program of
STM.\newline\indent \ \ E-mail: liucg@nankai.edu.cn
}~~~~~~~
~~~~~~~~~~
Youquan Zheng
\footnote
{School of Mathematical Sciences,  Nankai University,  Tianjin 300071,  China.\newline
 \indent \ \ E-mail: zhengyq@mail.nankai.edu.cn
}
}
\begin{document}
\date{}
\maketitle {\bf \noindent\large Abstract} In this paper, we consider a $p$-Laplacian equation in $\mathbf{R}^{N}$ with sign-changing potential and subcritical $p$-superlinear nonlinearity. By using the cohomological linking method for cones developed by Degiovanni and Lancelotti in 2007, an existence result is obtained. We also give a result on the existence of periodic solutions for one-dimensional $p$-Laplacian equations which can be proved by the same method.

\section{Introduction and main results}
We mainly consider the following $p$-Laplacian equation in the entire space
\begin{equation}\label{e1.1}
\left \{
\begin{array}
{ll} -\Delta_{p}u + U(x)|u|^{p-2}u= f(x, u), \\
u\in W^{1, p}(\mathbf{R}^{N}, \mathbf{R}),
\end{array}
\right.
\end{equation}
where $\Delta_{p}u:= div(|\nabla u|^{p-2}\nabla u)$ is the $p$-Laplacian operator with $ p >1$.

For $p=2$,  (\ref{e1.1}) turns into a kind of  Schr\"{o}dinger  equation of the form
\begin{equation}\label{e1.2}
- \Delta u + U(x)u = f(x, u), \;\;u\in H^1(\mathbf{R}^{N}, \mathbf{R}),
\end{equation}
which has been studied extensively.
In \cite{BPW01, BW95, BW00, LHW05}, (\ref{e1.2}) with  a constant sign potential $U(x)$ was considered. More precisely, the potential in these papers is of the form $a_{0}(x)+\lambda a(x)$, $a_{0}(x)$ has positive lower bound, $a(x)\geq 0$ and $\lambda > 0$ large enough. And in \cite{DS06, DS07, FT02}, the authors considered (\ref{e1.2}) with  a potential $U(x)$ that may change sign.

For general $p>1$, most of the work, as the authors of this paper known, deal  with the problem (\ref{e1.1}) with  a constant sign potential $U(x)$,  see for example \cite{F08, L10, KM08} and the reference therein.

In this paper,  we consider (\ref{e1.1}) with sign-changing potential and subcritical $p$-superlinear  nonlinearity,  moreover, periodic conditions on the potential and nonlinearity are not needed. Assume that $U(x)$ is of the form $b(x)-\lambda V(x)$,  here $\lambda$ is a real number,  $b(x)$,  $V(x)$ and $f(x, t)$ satisfy the following conditions
\begin{enumerate}
\item[{(B)}]  $b\in C(\mathbf{R}^{N}, \mathbf{R})$, $\displaystyle\inf_{x\in\mathbf{R}^{N}}b(x) \geq b_{0} > 0$,
$meas(\{x\in \mathbf{R}^{N}: b(x) \leq M\}) < \infty,\;\forall M>0$,
\item[{(V)}]   $V\in L^{\infty}(\mathbf{R}^{N}, \mathbf{R})$,
\item[{(f$_{1}$)}]   $f\in C(\mathbf{R}^{N}\times \mathbf{R}, \mathbf{R})$, $\exists \,q\in (p, p^{*})$ s.t. $|f(x, t)| \leq C(1 + |t|^{q-1})$,   $f(x, t)t \geq 0$, $\forall\,t\in\mathbf{R}$,
\item[{(f$_2$)}]  $\displaystyle\lim_{|t|\rightarrow \infty}\frac{f(x, t)t}{|t|^{p}} = + \infty$ uniformly in $x\in \mathbf{R}^{N}$,
\item[{(f$_3$)}]  $f(x, t) = o(|t|^{p-1})$ as $|t| \rightarrow 0$,  uniformly in $x\in \mathbf{R}^{N}$,
\item[ {(f$_4$)}] $\exists\theta \geq 1$ s.t. $\theta \mathcal{F}(x, t) \geq \mathcal{F}(x, st)$, $\forall\,(x, t)\in \mathbf{R}^{N} \times \mathbf{R}$ and $s\in [0, 1]$,
\end{enumerate}
here we have set $F(x, t)=\int_{0}^{t}f(x, t){\rm d}t$, $\mathcal{F}(x, t) = f(x, t)t - pF(x, t)$, $p^*=\frac{Np}{N-p}$ if $p<N$ and $p^*=+\infty$ if $p\ge N$ and $meas(\cdot)$ means the Lebesgue measure in $\mathbf{R}^{N}$.

Our main result reads as
\begin{thm}\label{t1.1}
If $(B)$, $(V)$ and  (f$_{1}$)-(f$_{4}$)  hold,  the problem $(\ref{e1.1})$ possesses a nontrivial solution for every $\lambda\in \mathbf{R}$.
\end{thm}

We remark that the condition $\displaystyle\inf_{x\in\mathbf{R}^{N}}b(x) \geq b_{0} > 0$ is not essential,  it can be replaced by the condition $\displaystyle\inf_{x\in \mathbf{R}^{N}}b(x) > -\infty$. First note that the case $\lambda = 0$ can be replaced by the case $V = 0$ with a nonzero $\lambda$,  we can always assume that $\lambda \neq 0$. If $\displaystyle\inf_{x\in\mathbf{R}^{N}}b(x)> -c_0$ for some $c_0>0$, one can replace $b$ and $V$ by $b + c_0$ and $V + \frac{c_0}{\lambda}$,   then $b + c_0$ and $V + \frac{c_0}{\lambda}$ satisfy conditions $(B)$ and $(V)$.
For $p=2$, the condition  (f$_4$)  was introduced in \cite{J99},  and for $p \neq 2$ it was introduced in \cite{LL03}.
Condition (B) was first introduced in \cite{BW95}, and then was used by many authors, for example, \cite{Z01}.

When dealing with superlinear problem, one usually needs a growth condition together with the following classical condition which was introduced by Ambrosetti and Rabinowitz in \cite{AR},
\begin{equation}\label{e1.3}
{\rm There \;exists\; }\mu > 2 \;{\rm such \;that} \;for \;u \neq 0\; {\rm and}\; x\in \mathbf{R}^{N}, 0 < \mu F(x, u)\leq uf(x, u). \end{equation}
Since then,  many authors tried to weaken this condition,  see \cite{DS07, F93, L10, LL03, LW04, LWZ06,SZ04}.
In \cite{LWZ06}  the authors obtained a weak solution of (\ref{e1.2}) under the following conditions
\begin{enumerate}
\item[(C$_{1}$)] $U(x)\in C(\mathbf{R}^{N}, \mathbf{R})$,  $\displaystyle\inf_{x\in\mathbf{R}^{N}}U(x)\geq U_{0} > 0$,  $U(x)$ is $1$-periodic in each
of $x_{i},\; i=1,\cdots,N$,
\item[(C$_{2}$)] $f(x, t)\in C^{1}$ is $1$-periodic in each of $x_{i},\; i=1,\cdots,N$,  $f'_{t}$ is a Caratheodory function and there exists $C > 0$,  such that
$|f'_{t}(x, t)|\leq C(1+|t|^{2^{*}-2})$,  $\displaystyle\lim_{|t|\rightarrow \infty}\frac{|f(x, t)|}{|t|^{2^{*}-1}}=0$,  uniformly in $x\in \mathbf{R}^{N}$,
\item[(C$_{3}$)] $f(x, t) = o(|t|)$,  as $|t|\rightarrow 0$,  uniformly in $x$,
\item[(C$_{4}$)] $\displaystyle\lim_{|u|\rightarrow \infty}\frac{F(x, u)}{u^{2}}=\infty$,  uniformly in $x$,
\item[(C$_{5}$)] $\frac{f(x, t)}{|t|}$ is strictly increasing in $t$.\end{enumerate}

\noindent And in \cite{L10}   the author  got a weak solution of (\ref{e1.1}) with the following assumptions
\begin{enumerate}
\item[(D$_{1}$)] $V\in C(\mathbf{R}^{N})$,  is $1$-periodic in $x_{i},\,i=1,\cdots,N$ and $0 < \alpha \leq V(x) \leq \beta < +\infty$,
\item[(D$_{2}$)] $f\in C(\mathbf{R}^{N} \times \mathbf{R})$ is $1$-periodic in $x_{i},\,i=1,\cdots,N$,  and
$\displaystyle\lim_{|t|\rightarrow\infty}\frac{f(x, t)}{|t|^{p^{*}-1}}=0$,
\item[(D$_{3}$)] $\displaystyle\lim_{|t|\rightarrow \infty}\frac{F(x, t)}{|t|^{p}}=+\infty$ uniformly in $x\in \mathbf{R}^{N}$,
\item[(D$_{4}$)] $f(x, t) = o(|t|^{p-2}t)$ as $|t|\rightarrow 0$,  uniformly in $x\in\mathbf{R}^{N}$,
\item[(D$_{5}$)] There exists $\theta \geq 1$ such that $\theta\mathcal{F}(x, t) \geq \mathcal{F}(x, st)$ for $(x, t)\in \mathbf{R}^{N}\times \mathbf{R}$ and $s\in [0, 1]$.\end{enumerate}

From above we can see (\ref{e1.3}) is weaken to (C$_{4}$) with the cost (C$_{5}$) and to (D$_{3}$) with the cost (D$_{5}$), respectively. And condition (D$_{5}$) is weaker than (C$_{5}$)(c.f.\cite{LL03}).  In our result, (f$_2$)takes place the condition (\ref{e1.3}) but we need (f$_4$).

We should also mention that there is another line to weaken (\ref{e1.3}). In \cite{DS07}  for   $\lambda$  large enough the authors  got a nontrivial solution of (\ref{e1.2}) with $U(x) = \lambda V(x)$ under the following conditions
 \begin{enumerate}
\item[(E$_{1}$)] $V\in C(\mathbf{R}^{N}, \mathbf{R})$,  $V$ is bounded below,  $V^{-1}(0)$ has nonempty interior,
    \item[(E$_{2}$)] their exists $M > 0$ such that the set $\{x\in\mathbf{R}^{N}| V(x) < M\}$ is nonempty and has finite measure,
\item[(E$_{3}$)] $f\in C(\mathbf{R}^{N}\times \mathbf{R}, \mathbf{R})$,  $F(x, u)\geq 0$ for all $(x, u)$,  $f(x, u)= o(u)$ uniformly in $x$ as $u \rightarrow 0$,
\item[(E$_{4}$)] $F(x, u)/u^{2} \rightarrow \infty$ uniformly in $x$ as $|u|\rightarrow \infty$,
\item[(E$_{5}$)] $\frac{1}{2}f(x, u)u - F(x, u) > 0$ whenever $u \neq 0$,
\item[(E$_{6}$)] $|f(x, u)|^{\tau} \leq a_{1}\Big(\frac{1}{2}f(x, u)u - F(x, u)\Big)|u|^{\tau}$ for some $a_{1} > 0$,  $\tau > \max\{1, N/2\}$ and all $(x, u)$ with $|u|$ large enough.\end{enumerate}
The authors of \cite{DS07} also proved that conditions (E$_{4}$)(E$_{5}$)(E$_{6}$)  are weaker than (\ref{e1.3}).

The main method used in the proof of Theorem \ref{t1.1} is the linking structure over cones which was developed in \cite{DL07}. We will use the Cerami condition
instead of $(PS)$ condition. This method  is also valid in finding periodic solutions for one-dimensional $p$-Laplacian equation. We will give  a brief  argument  in section 5 for this topic.

The  paper is organized as follows. In section 2,  we give the variational settings,  recall a critical point theorem and some important properties of cohomological index. In section 3,  an eigenvalue problem is studied. We get a divergent sequence of eigenvalues for this eigenvalue problem  by cohomological index theory. In section 4,  we give a proof of Theorem \ref{t1.1}. In section 5,  we state an existence result for the periodic solutions of one-dimensional $p$-Laplacian equation.
\section{Preliminaries}
Let ${\mathcal W}:= \{u\in W^{1, p}(\mathbf{R}^{N}, \mathbf{R}):\int_{\mathbf{R}^{N}}(|\nabla u|^{p} + b(x)|u|^{p}){\rm d}x < \infty\}$ with $b(x)$ satisfying the condition (B). Then ${\mathcal W}$ is a reflexive,  separable Banach space with norm $\|u\|= \left(\int_{\mathbf{R}^{N}}(|\nabla u|^{p} + b(x)|u|^{p}){\rm d}x\right)^{\frac{1}{p}}$.
From Gagliardo-Nirenberg inequality and H\"{o}lder inequality,  we have ${\mathcal W}\hookrightarrow L^{q}(\mathbf{R}^{N}, \mathbf{R})$ for $p \leq q \leq p^{*}$. Moreover,  we have the following compactness result. It was proved in $\cite{ZS06}$ in the case $p=2$. For the general case, the proof is similar. We give it here for reader's convenience.
\begin{lem}\label{l2.1}
${\mathcal W}\hookrightarrow\hookrightarrow L^{s}(\mathbf{R}^{N}, \mathbf{R})$ for $p \leq s < p^{*}$.
\end{lem}
\pf: Let $\{u_{n}\}\subset {\mathcal W}$ be a bounded sequence of ${\mathcal W}$ such that $u_{n}\rightharpoonup u$ weakly in ${\mathcal W}$. Then,  by the Sobolev embedding theorem,  $u_{n}\rightarrow u$ strongly in $L^{s}_{loc}(\mathbf{R}^{N}, \mathbf{R})$ for $p \leq s < p^{*}$. We first claim that
\begin{equation}\label{e2.3}u_{n} \rightarrow u\;\mbox{ strongly in}\; L^{p}(\mathbf{R}^{N}, \mathbf{R}).\end{equation}
In fact, by the uniformly convex properties of $L^{p}(\mathbf{R}^{N}, \mathbf{R})$,  we only need to prove that $\alpha_{n}:=\|u_{n}\|_{p} \rightarrow \|u\|_{p}$ (cf.$p_{295}$ in \cite{ding}). Assume,  up to subsequence,  that $\alpha_{n} \rightarrow \alpha$.

Set
$$\begin{array}{lll}
B_R = \{x\in \mathbf{R}^{N}:\, |x|<R\},\\
A(R, M) =  \{x\in \mathbf{R}^{N}\setminus B_{R}: b(x) \geq M\}, \\
B(R, M) =  \{x\in \mathbf{R}^{N}\setminus B_{R}: b(x) < M\}.
\end{array}$$
Then
\begin{eqnarray*}
\int_{A(R, M)}|u_{n}|^{p}{\rm d}x \leq \int_{\mathbf{R}^{N}}\frac{b(x)}{M}|u_{n}|^{p}{\rm d}x \leq \frac{\|u_{n}\|^{p}}{M}.
\end{eqnarray*}
Choose $t\in(1, \frac{p^{*}}{p})$ and $t'$ such that $\frac{1}{t} + \frac{1}{t'} = 1$. Then
\begin{eqnarray*}
\int_{B(R, M)}|u_{n}|^{p}{\rm d}x \leq \left(\int_{B(R, M)}|u_{n}|^{pt}{\rm d}x\right)^{\frac{1}{t}}(meas(B(R, M)))^{\frac{1}{t'}}
\leq C\|u_{n}\|^{p}(meas(B(R, M)))^{\frac{1}{t'}}.
\end{eqnarray*}
Since $\{\|u_{n}\|\}$ is bounded and condition $(B)$ holds,  we may choose $R$,  $M$ large enough such that
$\frac{\|u_{n}\|^{p}}{M}$ and $meas(B(R, M))$ are small enough. Hence, $\forall\varepsilon > 0$,  we have
\begin{eqnarray*}
\int_{\mathbf{R}^{N}\setminus B_{R}}|u_{n}|^{p}{\rm d}x = \int_{A(R, M)}|u_{n}|^{p}{\rm d}x + \int_{B(R, M)}|u_{n}|^{p}{\rm d}x < \varepsilon.
\end{eqnarray*}
Thus,
$$\begin{array}{ll}
\|u\|^{p}_{p}&= \|u\|^{p}_{L^{p}(B_{R})} + \|u\|^{p}_{L^{p}(\mathbf{R}^{N}\setminus B_{R})}\\& \geq \displaystyle\lim_{n \rightarrow \infty}\|u_{n}\|^{p}_{L^{p}(B_{R})}=\lim_{n\to \infty}(\|u_n\|^p-\|u_n\|^{p}_{L^{p}(\mathbf{R}^{N}\setminus B_{R})}) \geq \alpha^{p} - \varepsilon.
\end{array}$$

On the other hand,  let $\Omega$ be a arbitrary domain in $\mathbf{R}^{N}$,  then
\begin{eqnarray*}
\int_{\Omega}|u_{n}|^{p}{\rm d}x \leq \int_{\mathbf{R}^{N}}|u_{n}|^{p}{\rm d}x \rightarrow \alpha^{p},
\end{eqnarray*}
hence $\|u\|_{p} \leq \alpha$. Thanks to the arbitrariness of $\varepsilon$, we have
 $\alpha = \|u\|_{p}$. So (\ref{e2.3}) is proved.

Finally, it is easy to prove that $u_{n} \rightarrow u$ in $L^{s}(\mathbf{R}^{N}, \mathbf{R})$ for $p \leq s < p^{*}$. In fact, if $s\in (p,p^*)$, there is a number $\lambda\in (0,1)$ such that
$\frac 1s=\frac{\lambda}{p}+\frac{1-\lambda}{p^*}$. Then  by the H\"{o}lder inequality,
$$\|u_n-u\|_s^s=\int_{{\mathbf R}^n}|u_n-u|^{\lambda s}|u_n-u|^{(1-\lambda) s}{\rm d}x\le
\|u_n-u\|_p^{\lambda s}\|u_n-u\|_{p^*}^{(1-\lambda) s}.$$
Since $u_n$ is bounded in $L^{p^*}({\mathbf R}^n,\mathbf R)$ and $\|u_n-u\|_p\to 0$, we have $u_{n}\rightarrow u$ in $L^{s}(\mathbf{R}^{N}, \mathbf{R})$.
\qed

In the following,  we consider the $C^{1}$ functional $\Phi: {\mathcal W}\rightarrow \mathbf{R}$ defined by
\begin{equation}\label{e2.4}
\Phi(u) = \frac{1}{p}\int_{\mathbf{R}^{N}}(|\nabla u|^{p} +
b(x)|u|^{p}){\rm d}x -
\frac{\lambda}{p}\int_{\mathbf{R}^{N}}V(x)|u|^{p}{\rm d}x -
\int_{\mathbf{R}^{N}}F(x, u){\rm d}x.
\end{equation}
It is clear that critical points of $\Phi$ are weak
solutions of (\ref{e1.1}). In order to find a critical point of this functional, we will use the following critical point theorem.
It was proved in \cite{DL07},
where the functional was supposed to satisfy the $(PS)$ condition.
Recently,
in \cite{D09},
the author extended it to more general case
(the functional space is completely regular topological space or metric space).
If the functional space is a real Banach space,
according to the proof of Theorem 6.10 in \cite{D09},
the Cerami condition is sufficient for the compactness of the set of critical points at a fixed level
and the first deformation lemma to hold (see \cite{PAR}).
So this critical point theorem still hold under the Cerami condition.
\begin{lem}\label{l2.2}(\cite{DL07})
Let ${\mathcal W}$ be a real Banach space and let $C_{-}$,  $C_{+}$ be two
symmetric cones in ${\mathcal W}$ such that $C_{+}$ is closed in ${\mathcal W}$,  $C_{-}
\cap C_{+} = \{0\}$ and  $$i(C_{-}\setminus \{0\}) = i({\mathcal W}\setminus
C_{+}) = m < \infty.$$ Define the following four sets by
\begin{eqnarray*}
&&D_{-}=\{u\in C_{-}:\|u\| \leq r_{-}\}, \\
&&S_{+}=\{u\in C_{+}:\|u\|=r_{+}\}, \\
&&Q=\{u + te: u\in C_{-},  t \geq 0, \|u +te\| \leq r_{-}\}, \;\; e\in {\mathcal W}\setminus C_-, \\
&&H=\{u + te:u\in C_{-},  t \geq 0,  \|u + te\| = r_{-}\}.
\end{eqnarray*}
Then $(Q,  D_{-} \cup H)$ links $S_{+}$ cohomologically in dimension
$m+1$ over $\mathbf{Z_{2}}$. Moreover,  suppose $\Phi\in C^1({\mathcal W},
\mathbf{R})$ satisfying the Cerami condition,  and
$\displaystyle\sup_{x\in D_{-} \cup H}\Phi(x) <
\displaystyle\inf_{x\in S^{+}}\Phi(x)$,
$\displaystyle\sup_{x\in Q}\Phi(x) < \infty$. Then $\Phi$ has
a critical value $d \geq \displaystyle\inf_{x\in
S^{+}}\Phi(x)$.
\end{lem}

For convenience,  let us recall the definition and some properties
of the cohomological index of Fadell-Rabinowitz for a
$\mathbf{Z_{2}}$-set,  see \cite{FR77, FR78, PAR} for details. For
simplicity,  we only consider the usual $\mathbf{Z_{2}}$-action on a
linear space,  i.e.,  $\mathbf{Z_{2}}=\{1, -1\}$ and the action is the
usual multiplication. In this case,  the $\mathbf{Z_{2}}$-set $A$ is
a symmetric set with $-A=A$.

Let $E$ be a normed linear space. We denote by $\mathcal {S}(E)$ the
set of all symmetric subsets of  $E$  which do not contain the
origin of $E$. For $A\in \mathcal {S}(E)$,  denote $\bar{A} =
A/\mathbf{Z_{2}}$. Let $\rho: \bar{A} \rightarrow \mathbf{R}P^{\infty}$
be the classifying map and $\rho^{*}:
H^{*}(\mathbf{R}P^{\infty})=\mathbf{Z_{2}}[\omega] \rightarrow
H^{*}(\bar{A})$ the induced homomorphism of the cohomology rings.
The cohomological index of $A$,  denoted by $i(A)$,  is defined by
$\displaystyle\sup\{k \geq 1: \rho^{*}(\omega^{k-1}) \neq 0\}$.
We list some properties of the cohomological index here for further use  in this paper.
Let $A, B\in \mathcal {S}(E)$, there hold
\begin{enumerate}
\item[(i1)] ({\bf monotonicity}) if $A \subseteq B$,  then $i(A) \leq i(B)$,
\item[(i2)] ({\bf invariance}) if $h:A \rightarrow B$ is an odd homeomorphism,  then $i(A)=i(B)$,
\item[(i3)] ({\bf continuity}) if $C$ is a closed symmetric subset of $A$,  then
there exists a closed  symmetric neighborhood $N$ of $C$ in
$A$,  such that $i(N) = i(C)$,  hence the interior of $N$   in $A$ is also  a
neighborhood of $C$ in $A$ and
$i({\rm int} N)=i(C)$,
\item[(i4)] ({\bf neighborhood of zero}) if $V$ is bounded closed symmetric
neighborhood of the origin   in $E$,  then $i(\partial V) =
\dim E$.\end{enumerate}
\section{Eigenvalue problem}
In this section,  we consider the following eigenvalue problem
\begin{eqnarray*}
\left \{
\begin{array}
{ll} -\Delta_{p}u + b(x)|u|^{p-2}u=\lambda V(x)|u|^{p-2}u, \\
u\in W^{1, p}(\mathbf{R}^{N}, \mathbf{R}).\\
\end{array}
\right.
\end{eqnarray*}
We assume $V$ satisfies condition (V) and further assume that $V^{+}(x):=\frac{V(x)+|V(x)|}{2} \neq 0$ on some positive measure subset of ${\mathbf R}^N$ in this section.
Define on ${\mathcal W}$ the functionals
\begin{eqnarray*}
H(u)=\frac{1}{p}\int_{\mathbf{R}^{N}}(|\nabla u|^{p}+b(x)|u|^{p}){\rm d}x,
\end{eqnarray*}
\begin{eqnarray*}
I(u)=\frac{1}{p}\int_{\mathbf{R}^{N}}V(x)|u|^{p}{\rm d}x.
\end{eqnarray*}
Then
\begin{eqnarray*}
H\in C^{1}({\mathcal W}, \mathbf{R}),\;\;\;     \langle H'(u),  v\rangle = \int_{\mathbf{R}^{N}}(|\nabla u|^{p-2}\nabla u \nabla v
+ b(x)|u|^{p-2}uv){\rm d}x
\end{eqnarray*}
and
\begin{eqnarray*}
I\in C^{1}({\mathcal W}, \mathbf{R}), \;\;\;\;    \langle I'(u),  v\rangle = \int_{\mathbf{R}^{N}}V(x)|u|^{p-2}uv{\rm d}x.
\end{eqnarray*}
Our aim is to solve the eigenvalue problem
\begin{equation}\label{e3.3}
H'(u) = \lambda I'(u).
\end{equation}
\begin{lem}\label{l3.1} For any $u, v\in {\mathcal W}$,
it holds that \begin{eqnarray*}\langle
H'(u)-H'(v), u-v\rangle \geq
(\|u\|^{p-1}-\|v\|^{p-1})(\|u\|-\|v\|).\end{eqnarray*}
\end{lem}
\pf: We follow the idea of the proof of Lemma 2.3 in \cite{L06}. By direct computations,  we have
\begin{eqnarray*}
\langle H'(u)-H'(v), u-v\rangle &=&
\int_{\mathbf{R}^{N}}|\nabla u|^{p}+|\nabla v|^{p}-|\nabla u|^{p-2}\nabla
u\cdot \nabla v-|\nabla
v|^{p-2}\nabla v\cdot \nabla u{\rm d}x\\
& &+\int_{\mathbf{R}^{N}}b(x)
(|u|^{p}+|v|^{p}-|u|^{p-2}uv-|v|^{p-2}vu){\rm d}x.
\end{eqnarray*}
From the definition of the norm in ${\mathcal W}$,  we can get
\begin{eqnarray*} \langle
H'(u)-H'(v), u-v\rangle
&=& \|u\|^{p}+\|v\|^{p}
-\int_{\mathbf{R}^{N}}(|\nabla u|^{p-2}\nabla u\cdot \nabla v+b(x)
|u|^{p-2}uv){\rm d}x \\
& &-\int_{\mathbf{R}^{N}}(|\nabla v|^{p-2}\nabla v\cdot \nabla u+b(x)
|v|^{p-2}vu){\rm d}x.
\end{eqnarray*}
Applying H\"{o}lder inequality,
\begin{eqnarray*}
&&\int_{\mathbf{R}^{N}}(|\nabla u|^{p-2}\nabla u\cdot \nabla v+b(x)
|u|^{p-2}uv){\rm d}x
\\&&\leq\Big(\int_{\mathbf{R}^{N}}|\nabla
u|^{p}{\rm d}x\Big)^{\frac{p-1}{p}}\Big(\int_{\mathbf{R}^{N}}|\nabla
v|^{p}{\rm d}x\Big)^{\frac{1}{p}}+\Big(\int_{\mathbf{R}^{N}}b(x)
|u|^{p}{\rm d}x\Big)^{\frac{p-1}{p}}\Big(\int_{\mathbf{R}^{N}}b(x)
|v|^{p}{\rm d}x\Big)^{\frac{1}{p}}.
\end{eqnarray*}
Using the following
inequality
\begin{eqnarray*}(a+b)^{\alpha}(c+d)^{1-\alpha} \geq
a^{\alpha}c^{1-\alpha}+b^{\alpha}d^{1-\alpha}
\end{eqnarray*}
which holds for any $\alpha \in (0, 1)$ and for any $a > 0$,  $b > 0$,
$c
> 0$,  $d
> 0 $, set $\alpha = \frac{p-1}{p}$ and \begin{eqnarray*}
a=\int_{\mathbf{R}^{N}}|\nabla u|^{p}{\rm d}x,\;\;\;\;   b=\int_{\mathbf{R}^{N}}b(x)
|u|^{p}{\rm d}x,\;\;\;\;   c=\int_{\mathbf{R}^{N}}|\nabla v|^{p}{\rm d}x,\;\;\;\;
d=\int_{\mathbf{R}^{N}}b(x)|v|^{p}{\rm d}x,
\end{eqnarray*}
we can deduce that
\begin{eqnarray*} \int_{\mathbf{R}^{N}}(|\nabla
u|^{p-2}\nabla u\cdot \nabla v+b(x)|u|^{p-2}uv){\rm d}x \leq
\|u\|^{p-1}\|v\|.
\end{eqnarray*}
Similarly,  we can obtain
\begin{eqnarray*}
\int_{\mathbf{R}^{N}}(|\nabla v|^{p-2}\nabla v\cdot \nabla u+b(x)
|v|^{p-2}vu){\rm d}x \leq \|v\|^{p-1}\|u\|.
\end{eqnarray*}
Therefore,  we have  \begin{eqnarray*} \langle
H'(u)-H'(v), u-v\rangle &\geq &
\|u\|^{p} + \|v\|^{p} -
\|u\|^{p-1}\|v\| -
\|v\|^{p-1}\|u\|\\
&= & (\|u\|^{p-1} - \|v\|^{p-1})
(\|u\| - \|v\|).
\end{eqnarray*} \qed
\begin{lem}\label{l3.2} If $u_{n}\rightharpoonup u$ and  $\langle
H'(u_{n}), u_{n}-u\rangle \rightarrow 0$,  then
$u_{n}\rightarrow u$ in ${\mathcal W}$.
\end{lem}
\pf: Since ${\mathcal W}$ is a reflexive Banach
space, it is isometrically  isomorphic to a locally uniformly convex space, so as it was proved in \cite{ding}, weak convergence and norm convergence imply strong
convergence. Therefore we
only need to show that $\|u_{n}\| \rightarrow
\|u\|$.

Note that
\begin{eqnarray*}
\lim_{n\rightarrow \infty}\langle
H'(u_{n})-H'(u), u_{n}-u\rangle
=\lim_{n\rightarrow \infty}(\langle
H'(u_{n}), u_{n}-u\rangle - \langle
H'(u), u_{n}-u\rangle) = 0.
\end{eqnarray*}
By Lemma \ref{l3.1} we have
\begin{eqnarray*}
\langle H'(u_{n})-H'(u),  u_{n}-u\rangle
\geq
(\|u_n\|^{p-1}-\|u\|^{p-1})(\|u_{n}\|-\|u\|)
\geq 0.
\end{eqnarray*}
Hence $\|u_{n}\| \rightarrow \|u\|$ as
$n \rightarrow \infty$ and the assertion follows.\qed
\begin{lem}\label{l3.3} $I'$ is weak-to-strong continuous,  i.e.
$u_{n}\rightharpoonup u$ in ${\mathcal W}$ implies $I'(u_{n})\rightarrow
I'(u)$.
\end{lem}
\pf: This is a direct consequence of Theorem 1.22 in \cite{ZS06} and Lemma \ref{l2.1}.\qed
\begin{lem}\label{l3.4}
If $u_{n}\rightharpoonup u$ in ${\mathcal W}$,  then $I(u_{n}) \rightarrow I(u)$.
\end{lem}
\pf:
\begin{eqnarray*}
p|I(u_{n}) - I(u)| &=& |\langle
I'(u_{n}), u_{n}\rangle - \langle I'(u), u\rangle|\\ &=& |\langle
I'(u_{n}) - I'(u),  u_{n}\rangle + \langle I'(u),
u_{n} - u\rangle|\\
&\leq & \|I'(u_{n}) - I'(u)\| \|u_{n}\| + o(1).
\end{eqnarray*}
Because $u_{n} \rightharpoonup u$,  $u_{n}$ is bounded. From Lemma
\ref{l3.3},  we have $I(u_{n}) \rightarrow I(u)$.\qed

Set $\mathcal{M} = \{u\in {\mathcal W}: I(u) = 1\}$. Clearly,  $I(u)=\frac{1}{p}\langle
I'(u), u\rangle$,  so $1$ is a regular value of the functional $I$.
Hence by the implicit theorem, $\mathcal{M}$ is a $C^{1}$-Finsler manifold. It is complete,  symmetric,  since $I$ is continuous and
even. Moreover,  $0$ is not contained in $\mathcal{M}$,  so the trivial
$\mathbf{Z}_{2}$-action on $\mathcal{M}$ is free. Set
$\widetilde{H} = H|_{\mathcal{M}}$.
\begin{lem}\label{l3.5}
If $u\in\mathcal{M}$ satisfies $\widetilde{H}(u) =
\lambda$ and $\widetilde{H}'(u) = 0$,  then
$(\lambda, u)$ is a solution of (\ref{e3.3}).
\end{lem}
\pf: By Proposition 3.14.9 in \cite{PAR},  the norm of
$\widetilde{H}'(u)\in T^{*}_{u}\mathcal{M}$ is given
by $\|\widetilde{H}'(u)\|_{u}^{*} =
\displaystyle\displaystyle\min_{\mu\in \mathbf{R}}\|H'(u) - \mu
I'(u)\|^{*}$(here the norm $\|\cdot\|^{*}_{u}$ is the norm in the
fibre $T^{*}_{u}\mathcal{M}$,  and $\|\cdot\|^{*}$ is the operator
norm, the minimal can be attained was proved in Lemma 3.14.10 in \cite{PAR}). Hence there exists $\mu\in \mathbf{R}$ such that
$H'(u) - \mu I'(u) = 0$,  that is $(\mu,  u)$ is a
solution of (\ref{e3.3}) and $ \lambda=\widetilde{H}(u) = \frac{1}{p}\langle
H'(u), u\rangle = \frac{1}{p}\langle \mu I'(u), u\rangle = \mu \frac{1}{p}\langle
I'(u), u\rangle = \mu I(u)=
\mu$. \qed
\begin{lem}\label{l3.6} $\widetilde{H}$ satisfies the
$(PS)$ condition,  i.e. if $(u_{n})$ is a sequence on $\mathcal{M}$
such that $\widetilde{H}(u_{n}) \rightarrow c$,  and
$\widetilde{H}'(u_{n}) \rightarrow 0$,  then up to a
subsequence $u_{n} \rightarrow u\in \mathcal{M}$ in ${\mathcal W}$.
\end{lem}
\pf: First,  from the definition of $H$,  we can deduce
that $(u_{n})$ is bounded. Then,  up to a subsequence,
$u_{n}$ converges weakly to some $u$, by Lemma \ref{l3.4}, we have $I(u)=1$, so $u\in \mathcal{M}$.

From $\widetilde{H}'(u_{n}) \rightarrow 0$,  we have
$H'(u_{n}) - \mu_{n} I'(u_{n}) \rightarrow 0$ for a
sequence of real numbers $(\mu_{n})$. So  $\langle H'(u_{n}) - \mu_{n} I'(u_{n}), u_n\rangle \rightarrow 0$, thus we get $\mu_{n} \rightarrow c$. By Lemma \ref{l3.3},
we have $H'(u_{n})\rightarrow c I'(u)$. Hence $$\langle
H'(u_{n}),  u_{n} - u \rangle=\langle
H'(u_{n})-cI'(u),  u_{n} - u \rangle+\langle
cI'(u),  u_{n} - u \rangle \rightarrow 0.$$ By Lemma
\ref{l3.2},  we obtain $u_{n} \rightarrow u$. \qed

Let $\mathcal{F}$ denote the class of symmetric subsets of
$\mathcal{M}$, $\mathcal{F}_{n} = \{M\in \mathcal{F}: i(M) \geq
n\}$ and \begin{equation}\label{e3.4}\lambda_{n} = \displaystyle\inf_{M\in
\mathcal{F}_{n}}\displaystyle\sup_{u\in
M}\widetilde{H}(u).\end{equation} Since $\mathcal{F}_{n}
\displaystyle\supset \mathcal{F}_{n+1}$,  $\lambda_{n}
\leq \lambda_{n+1}$.
\begin{lem}\label{l3.7}For every $\mathcal{F}_{n}$, there is a  symmetric compact set $M\in \mathcal{F}_{n}$.
\end{lem}
\pf: We follow the idea of the proof of Theorem 3.2 in \cite{HT}. Since $meas\{x\in \mathbf{R}^{N}:V(x) > 0\} > 0$,
it implies  that $\forall n\in
\mathbf{N}^{*}$,  there exist $n$ open balls $(B_{i})_{1 \leq i \leq
n}$ in $\mathbf{R}^{N}$ such that $B_{i}\cap B_{j} = \emptyset$ for
$i \neq j$ and $meas(\{x\in
\mathbf{R}^{N}: V(x) > 0\}\cap B_{i}) > 0$. Approximating the characteristic
function $\chi_{\{x\in \mathbf{R}^{N}: V(x) > 0\}\cap B_{i}}$ by
$C^{\infty}(\mathbf{R}^{N}, \mathbf{R})$ functions in
$L^{p}(\mathbf{R}^{N}, \mathbf{R})$,  we can infer that there exists
a sequence $\{u_{i}\}_{1 \leq i \leq n}\subseteq
C^{\infty}(\mathbf{R}^{N}, \mathbf{R})$ such that
$\int_{\mathbf{R}^{N}}V(x)|u_{i}|^{p}{\rm d}x > 0 $ for all $i = 1, ...n$
and ${\rm supp}\,u_{i} \cap  {\rm supp}\,u_{j}=\emptyset$ when $i \neq j$. Normalizing $u_{i}$,  we assume that $I(u_{i}) = 1$. Denote
$U_{n}$ the space spanned by $(u_{i})_{1 \leq i \leq n}$. $\forall
u\in U_{n}$,  we have $u = \sum\limits_{i=1}^{n}\alpha_{i}u_{i}$ and
$I(u)=\sum\limits_{i=1}^{n}|\alpha_{i}|^{p}$. So $u\to
\Big(I(u)\Big)^{\frac{1}{p}}$ defines a norm on $U_{n}$. Since
$U_{n}$ is $n$ dimensional,  this norm is equivalent to
$\|\cdot\|$. Thus $\{u\in U_{n}: I(u)=1\}\subseteq
\mathcal{M}$ is compact with respect to the norm
$\|\cdot\|$ and by (i4),  $i(\{u\in
U_{n}:I(u)=1\})=n$. So $\{u\in U_{n}: I(u)=1\}\in \mathcal{F}_{n}$.\qed

By Lemma \ref{l3.7}, we have $\lambda_{n}<+\infty$, and by condition (B), there holds $\lambda_n\ge 0$. Furthermore,  by Lemma \ref{l3.6} and Proposition 3.14.7 in \cite{PAR},  we see that $\lambda_n$ is
sequence of critical values of $\widetilde{H}$ and  $\lambda_{n}\to +\infty,\; \mbox{as}\;n\to \infty$ . By
Lemma \ref{l3.5} we get a divergent sequence of eigenvalues for problem
$(\ref{e3.3})$. So we have the following result.
\begin{thm}\label{t3.8}
Problem $(\ref{e3.3})$ has an  increasing sequence eigenvalues $\lambda_n$ which are defined by (\ref{e3.4}) and  $\lambda_{n}\to +\infty,\; \mbox{as}\;n\to \infty$ .
\end{thm}
\begin{lem}\label{l3.9} Set
\begin{equation}\label{e3.5}\mu_{n} = \displaystyle\inf_{K\in \mathcal{F}^{c}_{n}}\displaystyle\sup_{u\in
K}H(u), \end{equation} where $\mathcal{F}^{c}_{n}=\{K\in \mathcal
{F}_n: K ~is~ compact \}$. Then we have $\lambda_n=\mu_n$.
\end{lem}
\pf: From Lemma \ref{l3.7}, $\mathcal{F}^{c}_{n}\ne \emptyset$ and so $\mu_n<+\infty$. It is obvious that $\lambda_n\le \mu_n$.  If $\lambda_n<\mu_n$, there is $M\in \mathcal{F}_{n}$ such that $\sup\limits_{u\in M}H(u)<\mu_n$. The closure $\overline{M}$ of $M$ in $\mathcal{M}$ is still in $\mathcal{F}_{n}$, by continuity of $H$, $\sup\limits_{u\in \overline{M}}H(u)<\mu_n$ holds. Applying the property  (i3) of the cohomological index, we can find a small open neighborhood  $A\in \mathcal{F}_{n}$ of $\overline{M}$ in $\mathcal{M}$ such that $\sup\limits_{u\in A}H(u)<\mu_n$. As it was proved in the proof of Proposition $3.1$ in
\cite{DL07},  for every symmetric   open subset $A$ of
$\mathcal{M}$,   there holds $i(A) = \displaystyle\sup\{i(K): K ~{\rm is ~compact~ and
~ symmetric~
\\with} ~ K \subseteq A\}$.  So we
can choose a symmetric compact subset $K\subseteq A$ with $i(K)\ge n$ and $\sup\limits_{u\in K}H(u)<\mu_n$. This contradicts to the definition of $\mu_n$. Therefore we have $\lambda_n=\mu_n$.\qed

Motivated by Theorem $3.2$ in \cite{DL07},  we have the following statement.
\begin{thm}\label{t3.10}If  $\lambda_{m} <
\lambda_{m+1}$ for some $m\in \mathbf{N}^{*}$,  then\\
$i(\{u\in {\mathcal W}\setminus \{0\}:H(u) \leq
\lambda_{m} I(u)\})$ = $i(\{u\in {\mathcal W}: H(u) <
\lambda_{m+1} I(u)\})$ = m.
\end{thm}
\pf: Suppose $\lambda_{m} < \lambda_{m+1}$.
If we set $A = \{u\in \mathcal{M}: H(u) \leq
\lambda_{m}\}$ and $B = \{u\in
\mathcal{M}:H(u) < \lambda_{m+1}\}$,
by the definition (\ref{e3.4}),  we have $i(A) \leq m$. Assume that $i(A) \leq m-1$. Thanks to (i3),
there exists a symmetric neighborhood $\mathcal{N}$ of $A$ in
$\mathcal{M}$ satisfying $i(\mathcal{N}) = i(A)$.
By the equivariant deformation theorem(see \cite{D03}),  there exists $\delta > 0$
and an odd continuous map $\iota:\{u\in \mathcal{M}:
H(u) \leq \lambda_{m}+\delta\} \rightarrow
\{u\in \mathcal{M}: H(u) \leq \lambda_{m}-\delta\}\cup \mathcal{N} = \mathcal{N}$. Hence $i(u\in \mathcal{M}:H(u) \leq
\lambda_{m}+\delta) \leq m-1$.
By (\ref{e3.4}), there exists $M\in \mathcal F_m$ such that $\sup\limits_{u\in M}H(u)<\lambda_m+\delta$. So $M\subseteq \{u\in \mathcal{M}:H(u) \leq
\lambda_{m}+\delta\}$ and thus $i(M)\le m-1$.
This contradicts to the fact that $M\in \mathcal F_m$. Thus we have $i(A) = m$. By the invariance of the cohomological index
under odd homeomorphism and the functionals $H, I$ are $p$-homogeneous,  we have $i(\{u\in {\mathcal W}\setminus
\{0\}:H(u) \leq \lambda_{m} I(u)\}) = m$.

Since $A\subseteq B$ and $i(A)=m$,  we have $i(B) \geq
m$. Assume that $i(B) \geq m+1$. As in the proof of Lemma \ref{l3.9},
there exists a symmetric,  compact subset $K$ of $B$ with $i(K) \geq
m+1$. Since $\max\limits_{ u\in K}H(u)<
\lambda_{m+1}=\mu_{m+1}$,  this contradicts to definition (\ref{e3.5}). By the
invariance of the cohomological index under odd homeomorphism and the functionals $H, I$ are $p$-homogeneous,  we
have $i(\{u\in {\mathcal W}: H(u) < \lambda_{m+1}
I(u)\}) = m $.\qed
\section{Proof of the main theorem}
Set $J(u) = \int_{\mathbf{R}^{N}}F(x, u){\rm d}x$,  by the definition of $H$ and $I$ in section 3,  we can write the functional $\Phi$ defined in section 2 as
\begin{eqnarray*}
\Phi(u) = H(u)-\lambda I(u)- J(u),\;\;u\in \mathcal W.
\end{eqnarray*}
It follows from Lemma 1.22 in \cite{ZS06} that
$J'$ is compact.

Replacing $(\lambda,  V)$ with $(-\lambda,  -V)$,  we can assume that $\lambda \geq 0$.

First,  we consider the case  $V^{+}(x) \neq 0$ on some positive measure subset of ${\mathbf R}^N$ and there exist $m \geq 1$ such that $\lambda_{m} \leq
\lambda <\lambda_{m+1}$. Set
\begin{equation}\label{e4.4}
C_{-} = \{u\in {\mathcal W}: H(u) \leq \lambda_{m}
I(u)\},\end{equation}
\begin{equation}\label{e4.5} C_{+} = \{u\in {\mathcal W}: H(u) \geq
\lambda_{m+1} I(u)\}.
\end{equation}
It is easy to see that $C_{-}$,  $C_{+}$ are two symmetric closed cones in ${\mathcal W}$
and $C_{-} \cap C_{+} = \{0\}$.
By Theorem \ref{t3.10} we have \begin{equation}\label{e4.1}i(C_{-}\setminus\{0\}) = i({\mathcal W}\setminus
C_{+}) = m.\end{equation}
\begin{thm}\label{t4.1}There exist $r_{+} > 0$ and $\alpha > 0$ such that $\Phi(u) > \alpha$ for $u\in C_{+}$ and $\|u\| = r_{+}$.\end{thm}
\pf: Let $\varepsilon > 0$ be small enough,  from (f$_1$) and (f$_3$),  we have $|F(x, t)| \leq \varepsilon|t|^{p}+C_{\varepsilon}|t|^{q}$,  by the Sobolev embedding inequality, for $u\in C_+$,  we can get
\begin{equation}\label{e4.2} \begin{array}{llllll}
\Phi(u)
&=& H(u)-\lambda I(u)- J(u) \\
&=& H(u) - \frac{\lambda}{\lambda_{m+1}}\lambda_{m+1}I(u)-J(u)
\\&\geq& H(u) - \frac{\lambda}{\lambda_{m+1}}H(u) - \varepsilon \int_{\mathbf{R}^{N}}|u|^{p}{\rm d}x -
C_{\varepsilon}\int_{\mathbf{R}^{N}}|u|^{q}{\rm d}x\\
&\geq& H(u) - \frac{\lambda}{\lambda_{m+1}}H(u) - \frac{\varepsilon}{b_{0}} \int_{\mathbf{R}^{N}}b(x)|u|^{p}{\rm d}x -
C_{\varepsilon}\int_{\mathbf{R}^{N}}|u|^{q}{\rm d}x\\
&\geq& (1 - \frac{\lambda}{\lambda_{m+1}} -  \frac{\varepsilon}{b_{0}})H(u)-
C_{\varepsilon}\int_{\mathbf{R}^{N}}|u|^{q}{\rm d}x\\
&\geq& \frac{(1 - \frac{\lambda}{\lambda_{m+1}} -  \frac{\varepsilon}{b_{0}})}{p}\|u\|^{p} - C\|u\|^{q}.
\end{array}\end{equation}
We remind that in the second inequality of (\ref{e4.2}), the condition (B) has been applied.
Since $p < q$,  the assertion follows.\qed

Since $\lambda \geq \lambda_{m}$,  by (f$_1$) it holds that
\begin{equation}\label{e4.3}\Phi(u) \leq 0,\;\forall\, u\in C_{-}.\end{equation}
 Set $\mathbf R^+=[0,+\infty)$. Following the idea of the proof of Theorem 4.1 in \cite{DL07}, we have
\begin{thm}\label{t4.2}
Let $e\in {\mathcal W}\setminus C_{-}$,  there exists $r_{-}
> r_{+}$ such that $\Phi(u) \leq 0$ for $u\in C_{-} +
\mathbf{R}^{+}e$ and $\|u\| \geq r_{-}$.\end{thm}
\pf: Define another norm on ${\mathcal W}$ by
$\|u\|_{V}:=(\int_{\mathbf{R}^{N}}(|V(x)|+1)|u|^{p}{\rm d}x)^{1/p}$.
Then the same reason as the proof of Theorem 4.1 in \cite{DL07}, there exists some constant $b > 0$ such that $\|u + te\| \leq b
\|u + te\|_{V}$ for every $u\in C_{-}$,  $t \geq 0$ and some $b > 0$. That is
\begin{equation}\label{e4.7}
\int_{\mathbf{R}^{N}}(|\nabla (u+te)|^{p} + b(x)|u+te|^{p}){\rm d}x \leq b^{p}\int_{\mathbf{R}^{N}}(|V(x)|+1)|u+te|^{p}{\rm d}x.
\end{equation}

Let $\{u_{k}\}$ be a sequence such that $\|u_{k}\|\rightarrow +\infty$ and $u_{k}\in
C_{-}+\mathbf{R}^{+}e$.
Set $v_{k}=\frac{u_{k}}{\|u_{k}\|}$,  then,  up to
a subsequence,  $\{v_{k}\}$ converges to some $v$ weakly in ${\mathcal W}$ and
a.e.in $\mathbf{R}^{N}$. Note that Lemma \ref{l3.4} is also true for functional $\int_{\mathbf{R}^{N}}(|V(x)|+1)|u|^{p}{\rm d}x,\;\;u\in\mathcal W$,  it follows from (\ref{e4.7}) that
$\int_{\mathbf{R}^{N}}(|V(x)|+1)|v|^{p}{\rm d}x \geq \frac{1}{b^{p}}$. So $|v| \neq 0$ on a  positive measure set $\Omega_0$. Since $\displaystyle\lim_{|t|\rightarrow\infty}\frac{f(x, t)t}{|t|^{p}} = +\infty$
implies $\displaystyle\lim_{|t|\rightarrow\infty}\frac{F(x, t)}{|t|^{p}} = +\infty$,  from (f$_2$) we have
\begin{eqnarray*}
\lim_{k \rightarrow
\infty}\frac{F(x, u_{k}(x))}{\|u_{k}\|^{p}} =
\lim_{k\rightarrow\infty}\frac{F(x, \|u_{k}\|v_{k}(x))}
{\|u_{k}\|^{p}|v_{k}(x)|^{p}} |v_{k}(x)|^{p}=+\infty, \;x\in \Omega_0.
\end{eqnarray*}
 By (f$_1$) and Fatou's lemma
we can get $$\frac{\int_{\mathbf{R}^{N}}F(x, u_{k})
{\rm d}x}{\|u_{k}\|^{p}}\rightarrow +\infty,  \;{\rm as}\;k\to \infty.$$  By the arbitrariness of the sequence $\{u_{k}\}$,
we have $\frac{\int_{\mathbf{R}^{N}}F(x, u)
{\rm d}x}{\|u\|^{p}}\rightarrow +\infty$ as $\|u\|\rightarrow +\infty$ and $u\in C_{-}+\mathbf{R}^{+}e$.
Noting that
$$\frac{\Phi(u)}{\|u\|^{p}} = \frac{1}{p} - \frac{\lambda I(u)}{\|u\|^{p}} - \frac{\int_{\mathbf{R}^{N}}F(x, u)
{\rm d}x}{\|u\|^{p}}$$
and by conditions (B) and (V),$$\left|\frac{I(u)}{\|u\|^{p}}\right| \leq \frac{C\int_{\mathbf{R}^{N}}|u|^{p}{\rm d}x}{\|u\|^{p}} \leq \frac{C\int_{\mathbf{R}^{N}}b(x)|u|^{p}{\rm d}x}{\|u\|^{p}}\leq C,$$
the assertion follows.
\qed

\begin{thm}\label{t4.3}$\Phi$ satisfies the Cerami condition. i.e., for any sequence $\{u_{k}\}$  in ${\mathcal W}$ satisfying
$(1 + \|u_{k}\|)\Phi'(u_{k}) \rightarrow 0$ and
$\Phi(u_{k})\rightarrow c$ possesses a convergent subsequence.\end{thm}
\pf: Let $\{u_{k}\}$ be a sequence in ${\mathcal W}$ satisfying
$(1 + \|u_{k}\|)\Phi'(u_{k}) \rightarrow 0$ and
$\Phi(u_{k})\rightarrow c$.
We claim that $\{u_{k}\}$ is bounded in ${\mathcal W}$. Otherwise, if $\|u_k\|\to \infty$,  we consider $w_{k}:=\frac{u_{k}}{\|u_{k}\|}$. Then,  up to subsequence,  we get $w_k \rightharpoonup w$ in ${\mathcal W}$,  $w_k\rightarrow w$ in $L^{s}(\mathbf{R}^{N})$ for $p \leq s < p^{*}$ and $w_k(x)\rightarrow w(x)$ a.e. $x\in \mathbf{R}^{N}$ as $k\to \infty$.
If $w \neq 0$ in ${\mathcal W}$,  since $\Phi'(u_{k})u_{k} \rightarrow 0$, that is to say
\begin{equation}\label{e4.8}
\int_{\mathbf{R}^{N}}(|\nabla u_{k}|^{p} + b(x)|u_{k}|^{p}){\rm d}x - \lambda\int_{\mathbf{R}^{N}}V(x)|u_{k}|^{p}{\rm d}x - \int_{\mathbf{R}^{N}}f(x, u_{k})u_{k}{\rm d}x \rightarrow 0,
\end{equation}
from condition (V), we have $\frac{|\int_{\mathbf{R}^{N}}V(x)|u_{k}|^{p}{\rm d}x|}{\|u_{k}\|^{p}} \leq C$,
so by dividing the left hand side of (\ref{e4.8}) with $\|u_k\|^p$ there holds \begin{equation}\label{e4.9}\left|\int_{\mathbf{R}^{N}}\frac{f(x, u_{k})u_{k}}{\|u_{k}\|^{p}}{\rm d}x\right| \leq C.\end{equation} On the other hand,  by Fatou's lemma and condition (f$_2$) we have
\begin{eqnarray*}
\int_{\mathbf{R}^{N}}\frac{f(x, u_{k})u_{k}}{\|u_{k}\|^{p}}{\rm d}x = \int_{\{w_{k} \neq 0\}}|w_{k}|^{p}\frac{f(x, u_{k})u_{k}}{|u_{k}|^{p}}{\rm d}x \rightarrow \infty,
\end{eqnarray*}
this contradicts to (\ref{e4.9}).

If $w=0$ in ${\mathcal W}$,  inspired by \cite{J99}, we choose $t_k\in [0,1]$ such that $\Phi(t_{k}u_{k}):= \displaystyle\max_{t\in [0, 1]}\Phi(tu_{k})$. For any $\beta> 0$ and $\tilde{w}_{k}:=(2p\beta)^{1/p}w_{k}$,    by Lemma \ref{l3.4} and the compactness of $J'$ we have that
\begin{eqnarray*}
\Phi(t_{k}u_{k}) \geq \Phi(\tilde{w}_{k}) = 2\beta - \lambda \int_{\mathbf{R}^{N}}V(x)|\tilde{w}_{k}|^{p}{\rm d}x - \int_{\mathbf{R}^{N}}F(x,  \tilde{w}_{k}){\rm d}x  \geq \beta,
\end{eqnarray*}
when $k$ is large enough, this implies that
\begin{equation}\label{e4.10}\lim\limits_{k \rightarrow \infty}\Phi(t_{k}u_{k}) = \infty.\end{equation} Since
$\Phi(0)=0,\;\Phi(u_k)\to c$,  we have $t_{k}\in (0, 1)$.
By the definition of $t_k$,
\begin{equation}\label{e4.11}\langle \Phi'(t_{k}u_{k}), t_{k}u_{k}\rangle = 0.\end{equation}
From (\ref{e4.10}), (\ref{e4.11}), we have
\begin{eqnarray*}
\Phi(t_{k}u_{k})-\frac 1p\langle \Phi'(t_{k}u_{k}), t_{k}u_{k}\rangle =\int_{\mathbf{R}^{N}}\left(\frac{1}{p}f(x, t_{k}u_{k})t_{k}u_{k} - F(x, t_{k}u_{k})\right){\rm d}x \rightarrow \infty.
\end{eqnarray*}
By  (f$_4$),  there exists $\theta\ge 1$ such that
\begin{equation}\label{e4.12}
\int_{\mathbf{R}^{N}}(\frac{1}{p}f(x, u_{k})u_{k} - F(x, u_{k})){\rm d}x
\geq \frac{1}{\theta}\int_{\mathbf{R}^{N}}(\frac{1}{p}f(x, t_{k}u_{k})t_{k}u_{k} - F(x, t_{k}u_{k})){\rm d}x \rightarrow \infty.
\end{equation}
On the other hand,
\begin{equation}\label{e4.13}
\int_{\mathbf{R}^{N}}(\frac{1}{p}f(x, u_{k})u_{k} - F(x, u_{k})){\rm d}x
= \Phi(u_{k}) - \frac{1}{p}\langle\Phi'(u_{k}),u_{k}\rangle\to c_0.
\end{equation}
(\ref{e4.12}) and (\ref{e4.13}) are contradiction. Hence $\{u_{k}\}$ is bounded in ${\mathcal W}$. So up to a subsequence, we can assume that $u_{k}\rightharpoonup u$ for some $\mathcal W$.

Since $\Phi'(u_{k}) = H'(u_{k})-
\lambda I'(u_{k}) - J'(u_{k})\to 0$ and $I'$, $J'$ are compact, we have that $H'(u_{k})\to \lambda I'(u)+ J'(u)$
 in ${\mathcal W}^{*}$. So $$\langle H'(u_k),u_k-u\rangle=\langle H'(u_k)-(\lambda I'(u)+ J'(u)),u_k-u\rangle+\langle \lambda I'(u)+ J'(u),u_k-u\rangle\to 0.$$
By Lemma \ref{l3.2}, $u_{k}\to u$ in $\mathcal W$. \qed

{\bf Proof of Theorem \ref{t1.1}} Define $D_{-}$,  $S_{+}$,  $Q$,  $H$ as lemma 2.2,  then
from Theorem \ref{t4.1}, $\Phi(u) \geq \alpha
> 0$ for every $u\in S_{+}$
 , from Theorem \ref{t4.2}, $\Phi(u) \leq 0$ for
every $u\in D_{-}\cup H$ and $\Phi$ is bounded on $Q$. Applying Theorem \ref{t4.3},  it follows
that $\Phi$ has a critical value $d \geq \alpha
> 0$. Hence $u$ is a nontrivial weak solution of (\ref{e1.1}).

For the cases  $0 \leq \lambda < \lambda_{1}$ or $V^{+}(x) \equiv 0$, set $C_{-} = \{0\}$ and $C_{+} = {\mathcal W}$,  it is easy to see that the arguments in this section are also valid. So we get a nontrivial solution and the
proof of Theorem \ref{t1.1} is complete.\qed
\section{Periodic problem for one-dimensional $p$-Laplacian equation}
In this section,  we state a result which can be proved by the same methods as in the proof of Theorem \ref{t1.1}. We only outline the main points. Our result reads as
\begin{thm}\label{t5.1}
If $p > 1$,  $V\in L^{\infty}(\mathbf{R}, \mathbf{R})$,  $f\in C(\mathbf{R}\times \mathbf{R}, \mathbf{R})$ satisfies (f$_{1}$)-(f$_{4}$),  both $V(t)$ and $f(t,u)$ are $1$-periodic in $t$,  then
\begin{equation}\label{e5.5}
\left \{
\begin{array}
{ll} -\Delta_{p}u + |u|^{p-2}u = \lambda V(t)|u|^{p-2}u + f(t, u), \\
u(0) = u(1),  u'(0) = u'(1).
\end{array}
\right.
\end{equation}
has a nontrivial solution for every $\lambda\in \mathbf{R}$.
\end{thm}

The periodic solution of $p$-Laplacian equation has been considered in
many papers,  for example,  \cite{BM07, BM08, MM00}. Up to the author's
knowledge,  Theorem \ref{t5.1} is  new.

Let ${\mathcal W}:= W^{1, p}(S^{1}, \mathbf{R})$ with the norm $\|u\| = (\int_{S^{1}}|\nabla u|^{p} + |u|^{p}{\rm d}t)^{\frac{1}{p}}$, here $S^1=\mathbf R/\mathbf Z$. Then ${\mathcal W}$ is a reflexive,  separable Banach space. And ${\mathcal W}$ can be embedded into $L^{q}(\mathbf{R}, \mathbf{R})$ for any $p \leq q < \infty$. As in section 3,  we consider the eigenvalue problem
\begin{eqnarray*}
\left \{
\begin{array}
{ll} -\Delta_{p}u + |u|^{p-2}u = \lambda V(t)|u|^{p-2}u , \\
u(0) = u(1),  u'(0) = u'(1).
\end{array}
\right.
\end{eqnarray*}
We can get a divergent sequence of eigenvalues defined by  $\lambda_{n} = \displaystyle\inf_{M\in
\mathcal{F}_{n}}\displaystyle\sup_{u\in
M}\int_{S^{1}}|\nabla u|^{p} + |u|^{p}{\rm d}t$ if $V^{+}(t) \neq 0$ on a positive measure subset of $S^1$, here $\mathcal{F}_{n}$ is the class of symmetrical
subsets with Fadell-Rabinowitz index greater than $n$ of $\mathcal{M}:=\{u\in {\mathcal W}: \int_{S^{1}}|u|^{p}{\rm d}t = 1\}$. And if $\lambda_{m} <
\lambda_{m+1}$ for some $m\in \mathbf{N}^{*}$,  then
$i(\{u\in {\mathcal W}\setminus \{0\}: \int_{S^{1}}|\nabla u|^{p} + |u|^{p}{\rm d}t \leq
\lambda_{m} \int_{S^{1}}|u|^{p}{\rm d}t\})$ = $i(\{u\in {\mathcal W}: \int_{S^{1}}|\nabla u|^{p} + |u|^{p}{\rm d}t <
\lambda_{m+1} \int_{S^{1}}|u|^{p}{\rm d}t\}) = m$. Then arguing as in section 4,
consider the functional on ${\mathcal W}$
\begin{eqnarray*}
\Phi(u) = \frac{1}{p}\int_{S^{1}}(|\nabla u|^{p} + |u|^{p}){\rm d}t -
\frac{\lambda}{p}\int_{S^{1}}V(t)|u|^{p}{\rm d}t -
\int_{S^{1}}F(t, u){\rm d}t,
\end{eqnarray*}
assume $\lambda \geq 0$,  $V^{+}(t) \neq 0$ on a positive measure subset of $S^1$ and there exists $m\in \mathbf{N}^{*}$ such that $\lambda_{m} \leq \lambda < \lambda_{m+1}$,
Set
\begin{eqnarray*}
C_{-} = \{u\in {\mathcal W}: \int_{S^{1}}|\nabla u|^{p} + |u|^{p}{\rm d}t\leq \lambda_{m}
\int_{S^{1}}|u|^{p}{\rm d}t\},
\end{eqnarray*}
\begin{eqnarray*}
C_{+} = \{u\in {\mathcal W}:  \int_{S^{1}}|\nabla u|^{p} + |u|^{p}{\rm d}t\geq
\lambda_{m+1} \int_{S^{1}}|u|^{p}{\rm d}t\},
\end{eqnarray*}
then we have
\begin{enumerate}
\item[(1)]
There exist $r_{+} > 0$ and $\alpha > 0$ such that $\Phi(u) > \alpha$ for $u\in C_{+}$ and $\|u\| = r_{+}$,
\item[(2)] Let $e\in {\mathcal W}\setminus C_{-}$,  there exists $r_{-}
> r_{+}$ such that $\Phi(u) \leq 0$ for $u\in C_{-} +
\mathbf{R}^{+}e$ and $\|u\| \geq r_{-}$,
\item[(3)] $\Phi$ satisfies the Cerami condition.
\end{enumerate}
Then from Lemma \ref{l2.2},  we can get a nontrivial solution for (\ref{e5.5}). The cases for $0 \leq \lambda < \lambda_{1}$ or $V^{+}(x) \equiv 0$ are similar as in the proof of Theorem \ref{t1.1}.

\end{document}